# LINEAR, SECOND ORDER AND UNCONDITIONALLY ENERGY STABLE SCHEMES FOR THE VISCOUS CAHN-HILLIARD EQUATION WITH HYPERBOLIC RELAXATION USING THE INVARIANT ENERGY QUADRATIZATION METHOD

XIAOFENG YANG [†] AND JIA ZHAO [‡]

ABSTRACT. In this paper, we consider numerical approximations for the viscous Cahn-Hilliard equation with hyperbolic relaxation. This type of equations processes energy-dissipative structure. The main challenge in solving such a diffusive system numerically is how to develop high order temporal discretization for the hyperbolic and nonlinear terms, allowing large time-marching step, while preserving the energy stability, i.e. the energy dissipative structure at the time-discrete level. We resolve this issue by developing two second-order time-marching schemes using the recently developed "Invariant Energy Quadratization" approach where all nonlinear terms are discretized semi-explicitly. In each time step, one only needs to solve a symmetric positive definite (SPD) linear system. All the proposed schemes are rigorously proven to be unconditionally energy stable, and the second-order convergence in time has been verified by time step refinement tests numerically. Various 2D and 3D numerical simulations are presented to demonstrate the stability, accuracy and efficiency of the proposed schemes.

## 1. Introduction

The classical Cahn-Hilliard (CH) equation dates back to 1958 in Cahn and Hillard's seminal paper [5]. In the past decades, it has been well studied and broadly used to investigate the coarsening dynamics of two immersible fluids. Recently, researchers have devoted tremendous attention on the relaxed CH system, i.e. the viscous Cahn-Hilliard (VCH) system and its perturbed form with the hyperbolic relaxation (HR) effect (referred as to the perturbed viscous Cahn-Hilliard equation). Both VCH and VCH-HR have been well-studied theoretically where the topics are mainly focused on the well-posedness, sharp interface limit or global attractor, etc., see [2, 3, 6, 11–13, 16, 17, 19–21, 25–32, 36, 38, 43, 45–47, 56] and the references therein. Formally, the governing equation of the VCH-HR system is slightly different from the CH equation by incorporating two extra terms, including a strong damping (or called "viscosity") term and a hyperbolic relaxation term (or called "inertia"). The viscous effect is first proposed by Novick-Cohen [45] in order to introduce an additional regularity and some parabolic smoothing effects, can be viewed as a singular limit of the phase field equations for phase transitions [20]. The hyperbolic relaxation term was proposed by Galenko et. al. in [25–30, 38], in order to describe strongly non-equilibrium decomposition generated by rapid solidification under supercooling into the spinodal region occurring in certain materials (e.g., glasses). Since the VCH-HR system combines the hyperbolic relaxation and the viscosity together, it is mathematically more tractable comparing to the CH or VCH systems [36, 47, 78].

Before developing efficient numerical schemes to solve the VCH-HR system, we notice that its reduced version, the classical CH equation is now widely applied to model the interfacial dynamics in various scientific fields (cf. [5, 7, 18, 37, 39, 42, 57, 71, 72] and the references therein). The CH equation as well as its analogous counterpart model, the Allen-Cahn equation, are both categorized as representative equations of phase field type models. From the numerical point of view, when solving phase field models, it is desirable to establish efficient numerical schemes that can verify the so called "energy stable" property at the discrete level irrespectively of the coarseness of the discretization. In what follows, those algorithms will be called *unconditionally energy stable* or *thermodynamically consistent*. Schemes with this property are especially preferred since it is not only critical for the numerical scheme to capture the correct long time dynamics of the system, but also provides sufficient flexibility for dealing with the stiffness issue. In spite of this, we have to point out a basic fact that larger time step will definitely induce larger computational errors. In other words, the schemes with unconditional energy stability can allow arbitrary large time step only for the sake of the stability concern. In practice, the controllable accuracy is one of the most important factors to measure whether a scheme is reliable or not. Therefore, if one attempts to use the time step as large as possible while maintaining desirable accuracy, the only possible choice to develop more accurate schemes, e.g., the unconditionally energy stable second order schemes, which is the main focus of this paper.

It is remarkable that, despite a great deal of work done for the numerical solution of the classical CH system, almost all research related to the VCH or VCH-HR system had been focused on the theoretical PDE analysis with very few numerical analysis or algorithm design. More precisely, to the best of the authors' knowledge, no schemes can be claimed to posses the following three properties, namely, easy-to-implement, unconditionally energy stability and second order accuracy for the VCH-HR model. This is due to the numerical difficulties of proper discretization for the viscous effect and the hyperbolic inertia, besides for the regular stiffness issue induced by the nonlinear double well potential. At the very least, even for the reduced version, i.e. the CH system, the algorithm design is still challenging. It can be seen clearly from the following fact that some severe stability restrictions on the time step will occur if the nonlinear term is discretized in some normal ways like fully explicit type approach. Such a time step constraint can cause very high computational cost in practice [4, 23, 51]. Many efforts (primarily for CH system) had been done in order to remove this constraint and two commonly used techniques were developed, namely, the nonlinear convex splitting approach [22, 33, 56, 58], and the linear stabilized approach [8, 40, 41, 44, 48–55, 61, 62, 70, 73, 76, 77]. The convex splitting approach is unconditionally energy stable, but it produces nonlinear schemes, thus the implementation is complicated and the computational cost might be high. The linear stabilized approach is linear so it is efficient and very easy to implement. But, its stability requests a special property (generalized maximum principle) satisfied by the classical PDE solution or the numerical solution, which is not trivial to prove. Moreover, it is difficult to extend to second-order while preserving unconditional energy stability (cf. [51]).

Therefore, in order to develop some more efficient and accurate time marching schemes for solving the VCH-HR equation, we use the *Invariant Energy Quadratization* (IEQ) approach, which has been successfully applied to solve a variety of phase field type models, see [9, 10, 35, 59, 60, 63–66, 68, 69, 74, 75]). Its idea is very simple but quite different from those traditional methods like implicit, explicit, nonlinear splitting, or other various tricky Taylor expansions to discretize the nonlinear potentials. The essential strategy of IEQ is to make the free energy *quadratic*. To be more specific, the free energy potential is transformed into the quadratic form forcefully via the change of variables. Then, upon a simple reformulation, all nonlinear terms are treated by the semi-explicit way, which in turn yields a linear system. We develop two second order schemes, in which, one is based on the Crank-Nicolson, and the other is based on the Adam-Bashforth (BDF2). The schemes are *second order accurate*, *easy-to-implement* (linear system), and *unconditionally energy stable* (with a discrete energy dissipation law). Moreover, we show that the linear operator of all schemes are *symmetric positive definite*, so that one can solve it using the well-developed fast matrix solvers efficiently (CG or other Krylov subspace methods). Through various 2D and 3D numerical simulations, we demonstrate stability and accuracy of the proposed schemes.

The rest of the paper is organized as follows. In Section 2, we present the whole system and show the energy law in the continuous level. In Section 3, we develop the numerical schemes and prove their unconditional energy stabilities. In Section 4, we present various 2D and 3D numerical experiments to validate the accuracy and efficiency of the proposed numerical schemes. Finally, some concluding remarks are presented in Section 5.

## 2. Model Equations

First of all, we give a brief description for the model equations. We consider a binary alloy in a bounded domain $\Omega \in \mathbb{R}^d, d = 2, 3$ with $\partial \Omega$ Lipschitz continuous. For any $g_1, g_2 \in L^2(\Omega)$, we denote the inner product and $L^2$ norm as

$$(2.1) \quad (g_1, g_2) = \int_\Omega g_1 g_2 d\boldsymbol{x}, \quad \|g_1\| = \int_\Omega |g_1|^2 d\boldsymbol{x}.$$

We define $\phi(\boldsymbol{x}, t)$ as volume fraction of one material component, and $\boldsymbol{J}$ the diffusion flux, then the balance law for volume fraction gives

$$(2.2) \quad \phi_t + \nabla \cdot \boldsymbol{J} = 0.$$

In order to describe the evolution for $\phi$, we need to introduce a constitutive assumption on $\boldsymbol{J}$. One case could be

$$(2.3) \quad \alpha \boldsymbol{J}_t + \boldsymbol{J} = -\nabla(\frac{\delta E}{\delta \phi} + \beta \phi_t),$$

where $\alpha \geq 0$ is the relaxation parameter, $\beta \geq 0$ is the viscosity parameter, and $E(\phi)$ is the total free energy that takes the form as

$$(2.4) \quad E(\phi) = \int_\Omega \Big(\frac{\epsilon^2}{2}|\nabla \phi|^2 + F(\phi)\Big) d\boldsymbol{x},$$

where $\epsilon$ is a the positive constant that measures the interfacial width and $F(\phi)$ is a nonlinear bulk potential.

For the choice of nonlinear potential $F(\phi)$, we can choose either (i) double well (Ginzburg-Landau) potential where

$$(2.5) \quad F(\phi) = \phi^2(\phi - 1)^2;$$

or (ii) Flory-Huggins potential (cf. [47]) where

$$(2.6) \quad F(\phi) = (1-\phi)\ln(1-\phi) + \phi\ln\phi + \theta\phi(1-\phi), \quad \theta > 0.$$

By combining (2.2) and (2.3), the governing PDE reads as follows,

$$(2.7) \quad \alpha \phi_{tt} + \phi_t = \lambda \Delta \mu,$$
$$(2.8) \quad \mu = -\epsilon^2 \Delta \phi + f(\phi) + \beta \phi_t,$$

where $f(\phi) = F'(\phi)$, i.e., $f(\phi) = 2\phi(\phi-1)(2\phi-1)$ for double well potential and $f(\phi) = \ln(\frac{\phi}{1-\phi}) + \theta(1-2\phi)$ for Flory-Huggins potential. The boundary conditions can be

$$(2.9) \qquad (i) \text{ all variables are periodic; or } (ii) \ \partial_{\mathbf{n}}\phi|_{\partial\Omega} = \partial_{\mathbf{n}}\mu|_{\partial\Omega} = 0,$$

where $\mathbf{n}$ is the unit outward normal on $\partial\Omega$.

When $\alpha = \beta = 0$, the system degenerates to the standard CH system that conserves the local mass density. When $\alpha \neq 0$, the volume conservation will only hold provided $\int_\Omega \phi_t(0,\bm{x})d\bm{x} = 0$. To see this, by taking the $L^2$ inner product of (2.7) with 1, one can obtain directly

$$(2.10) \qquad \alpha \frac{d}{dt}\int_\Omega \phi_t(t,\bm{x})d\bm{x} + \int_\Omega \phi_t(t,\bm{x})d\bm{x} = 0.$$

This actually is an ODE system for time, and its solution is

$$(2.11) \qquad \int_\Omega \phi_t(t,\bm{x})d\bm{x} = \exp(\frac{-t}{\alpha})\int_\Omega \phi_t(0,\bm{x})d\bm{x}.$$

Therefore, by setting $\int_\Omega \phi_t(0,\bm{x})d\bm{x} = 0$, we obtain

$$(2.12) \qquad \int_\Omega \phi_t(t,\bm{x})d\bm{x} = \int_\Omega \phi_{tt}(t,\bm{x})d\bm{x} = 0.$$

Define the inverse Laplace operator $\Delta^{-1}$ such that $v = \Delta^{-1}u$ (with $\int_\Omega u d\bm{x} = 0$), iff

$$(2.13) \qquad \begin{cases} \Delta v = u, \quad \int_\Omega u d\bm{x} = 0, \\ \text{with the boundary conditions either (i) } v \text{ is periodic, or (ii) } \partial_{\mathbf{n}}v|_{\partial\Omega} = 0. \end{cases}$$

We now derive that the energy dissipation law for model (2.7)-(2.8). Denote the total energy

$$(2.14) \qquad \mathcal{E} = \int_\Omega \Big(\frac{\epsilon^2}{2}|\nabla\phi|^2 + F(\phi) + \frac{\alpha}{2}|\nabla\Delta^{-1}\phi_t|^2\Big)d\bm{x}.$$

Then, the model (2.7)-(2.8) satisfies the following energy dissipation law

$$(2.15) \qquad \frac{d\mathcal{E}}{dt} = -\|\nabla\Delta^{-1}\phi_t\|^2 - \beta\|\phi_t\|^2 \leq 0.$$

We introduce a new variable $\psi = \phi_t$. Since $\int_\Omega \psi d\bm{x} = \int_\Omega \psi_t d\bm{x} = 0$, using the operator $\Delta^{-1}$, we rewrite the system (2.7)-(2.8) as follows,

$$(2.16) \qquad \alpha\Delta^{-1}\psi_t + \Delta^{-1}\psi = -\epsilon^2\Delta\phi + f(\phi) + \beta\phi_t.$$

By taking the $L^2$ inner product of (2.16) with $\phi_t$, we have

$$(2.17) \qquad \alpha(\Delta^{-1}\psi_t, \psi) + (\Delta^{-1}\psi, \psi) - \beta\|\phi_t\|^2 = \frac{d}{dt}\int_\Omega \Big(\frac{\epsilon^2}{2}|\nabla\phi|^2 + F(\phi)\Big)d\bm{x}.$$

Since $\int_\Omega \psi d\bm{x} = 0$, we can find another auxillary variable $p$ such that $p = \Delta^{-1}\psi$, i.e.,

$$(2.18) \qquad \Delta p = \psi, \quad \int_\Omega \psi d\bm{x} = 0,$$

with the boundary condition specified in (2.13). By taking the $L^2$ inner product of (2.18) with $p$, that is

$$(2.19) \qquad (p, \psi) = -\|\nabla p\|^2 = (\Delta^{-1}\psi, \psi).$$

We differentiate (2.18) with time $t$ to obtain

$$(2.20) \qquad \Delta p_t = \psi_t.$$

By taking the $L^2$ inner product of (2.20) with $p$, we obtain $(\psi_t, p) = (\Delta p_t, p) = -\frac{1}{2} d_t \|\nabla p\|^2$. Hence, we derive

$$(2.21) \qquad \alpha(\Delta^{-1}\psi_t, \psi) = \alpha(\psi_t, \Delta^{-1}\psi) = \alpha(\psi_t, p) = -\frac{\alpha}{2} d_t \|\nabla p\|^2.$$

By combining (2.17)-(2.19)-(2.21), we obtain

$$(2.22) \qquad \frac{d}{dt} \int_\Omega \Big(\frac{\epsilon^2}{2}|\nabla \phi|^2 + F(\phi) + \frac{\alpha}{2}|\nabla p|^2\Big) d\bm{x} = -\|\nabla p\|^2 - \beta\|\phi_t\|^2 \leq 0,$$

that means the total free energy of the VCH-HR system (2.7)-(2.8) decays in time.

## 3. Numerical Schemes

We now construct two semi-discrete time marching numerical schemes for solving the model system (2.7)-(2.8)-(2.9) and prove their energy stabilities based on the *Invariant Energy Quadratization* (IEQ) approach. The intrinsic idea of the IEQ method is to transform the nonlinear potential into quadratic form. It is feasible since we notice that the nonlinear potential $F(\phi)$ is always bounded from below, in either the double well form (for the Ginzberg-Landau potential) or logarithmic form (for the Flory-Huggins potential). Thus, in general, we could rewrite the free energy functional $F(\phi)$ into the following equivalent form

$$(3.1) \qquad F(\phi) = (F(\phi) + B) - B,$$

where B is some constant to ensure $F(x) + B > 0, \forall x \in \mathbb{R}$, and define an auxilliary function $U$ as

$$(3.2) \qquad U = \sqrt{F(\phi) + B}.$$

Thus the total energy of (2.4) turns into a new form

$$(3.3) \qquad E(\phi, U) = \int_\Omega \Big(\frac{\epsilon^2}{2}|\nabla \phi|^2 + U^2 - B\Big) d\bm{x}.$$

Then we obtain an equivalent PDE system by taking the time derivative for the new variable $U$:

$$(3.4) \qquad \alpha\psi_t + \psi = \Delta\mu,$$
$$(3.5) \qquad \mu = -\epsilon^2 \Delta\phi + UH + \beta\phi_t,$$
$$(3.6) \qquad U_t = \frac{1}{2}H\phi_t,$$
$$(3.7) \qquad \psi = \phi_t,$$

where

$$(3.8) \qquad H(\phi) = \frac{f(\phi)}{\sqrt{F(\phi) + B}}, \; f(\phi) = F'(\phi).$$

The boundary conditions for the new system are still (2.9) since the equation (3.6) for the new variable $U$ is simply an ODE with time. The initial conditions read as

$$(3.9) \qquad \phi|_{(t=0)} = \phi_0, \psi|_{(t=0)} = 0,$$
$$(3.10) \qquad U|_{(t=0)} = \sqrt{F(\phi_0) + B},$$

where we simply set the initial profile of $\psi$ to be zero point-wise.

It is clear that the new transformed system (3.4)-(3.7) still retains a similar energy dissipative law. By applying the inverse Laplace operator $\Delta^{-1}$ to (3.4), taking the $L^2$ inner product of it with $\phi_t$, of (3.6) with $-U$, using (2.19) and (2.21), and summing them up, we can obtain the energy dissipation law of the new system (3.4)-(3.5). If we denote the transformed (equivalent) energy as

$$\hat{\mathcal{E}} = \int_\Omega \Big(\frac{\epsilon^2}{2}|\nabla\phi|^2 + U^2 + \frac{\alpha}{2}|\nabla\Delta^{-1}\psi|^2\Big)d\boldsymbol{x} \tag{3.11}$$

the energy law of (3.4)-(3.5) reads

$$\frac{d}{dt} = -\|\nabla\Delta^{-1}\phi_t\|^2 - \beta\|\psi\|^2 \leq 0. \tag{3.12}$$

**Remark 3.1.** *We emphasize that the new transformed system (3.4)-(3.7) is exactly equivalent to the original system (2.7)-(2.8), since (3.2) can be easily obtained by integrating (3.6) with respect to the time. For the time-continuous case, the potentials in the new free energy (3.3) are the same as the Lyapunov functional in the original free energy of (2.4), and the new energy law (3.12) for the transformed system is also the same as the energy law (2.15) for the original system as well. We will develop unconditionally energy stable numerical schemes for time stepping of the transformed system (3.4)-(3.7), and the proposed schemes should formally follow the new energy dissipation law (3.12) in the discrete sense, instead of the energy law for the originated system (2.15).*

**Remark 3.2.** *If $F(\phi) = \phi^2(\phi-1)^2$, we let $B = 0$, thus $H(\phi) = 2\phi-1$. At this time, the IEQ method is exactly the same as the so-called Lagrange multiplier method developed in [34]. We remark that the Lagrange multiplier method in [34] only works for the fourth order polynomial potential ($\phi^4$). This is because the term $\phi^3$ (the first order derivative of $\phi^4$) can be decomposed into a multiplication of two factors: $\lambda(\phi)\phi$, where $\lambda(\phi) = \phi^2$. In [34], this Lagrange multiplier term $\lambda(\phi)$ is then defined as the new auxiliary variable $U$. However, for other type potentials, e.g., the F-H potential, the new variable $U$ will takes a form as $\lambda(\phi) = \frac{1}{\phi}\ln(\frac{\phi}{1-\phi})$, this is unworkable for algorithms design.*

**Remark 3.3.** *If $F(\phi) = (1-\phi)\log(1-\phi) + \phi\log\phi + \theta\phi(1-\phi)$, following the work in [14], we regularize the logarithmic bulk potential by a $C^2$ piecewise function. More precisely, for any $0 < \sigma \ll 1$, the regularized free energy is*

$$\widehat{F}(\phi) = \begin{cases} \phi\ln\phi + \frac{(1-\phi)^2}{2\sigma} + (1-\phi)\ln\sigma - \frac{\sigma}{2} + \theta\phi(1-\phi), & \text{if } \phi \geq 1-\sigma, \\ \phi\ln\phi + (1-\phi)\ln(1-\phi) + \theta\phi(1-\phi), & \text{if } \sigma \leq \phi \leq 1-\sigma, \\ (1-\phi)\ln(1-\phi) + \frac{\phi^2}{2\sigma} + \phi\ln\sigma - \frac{\sigma}{2} + \theta\phi(1-\phi), & \text{if } \phi \leq \sigma. \end{cases} \tag{3.13}$$

*For convenience, we consider the problem formulated with the substitute $\widehat{F}(\phi)$, but omit the $\widehat{\phantom{x}}$ in the notation. Now the regularized functional $F(\phi)$ is defined in $\mathbb{R}$. Small fluctuation of the numerical solution $\phi$ near the boundary $(0,1)$ would not cause blow up of the numerical solution. In [14], the authors proved the error bound between the regularized PDE and the original PDE is controlled by $\sigma$ up to a constant. For this case, we simply take $B = 1$ that can ensure $\widehat{F}(x) + B > 0, \forall x \in \mathbb{R}$.*

The time marching numerical schemes are developed to solve the new transformed system (3.4)-(3.7). The proof of the unconditional stability of the schemes follows the similar lines as in the derivation of the energy law (3.12). Let $\delta t > 0$ denote the time step size and set $t^n = n\,\delta t$ for $0 \leq n \leq N$ with the ending time $T = N\,\delta t$.

### 3.1. Crank-Nicolson Scheme.
We first develop a second order scheme based on Crank-Nicolson method, that reads as follows.

<u>Scheme</u> **1.** *Given the initial condition $(U^0, \phi^0)$, compute $U^1$ and $\phi^1$ by assuming $U^{-1} = U^0$ and $\phi^{-1} = \phi^0$ for the initial step. Having computed $(\phi^n, U^n)$ and $(\phi^{n-1}, U^{n-1})$, with $n \geq 1$, we update*

$\phi^{n+1}$ and $U^{n+1}$ as follows:

$$\alpha \frac{\psi^{n+1} - \psi^n}{\delta t} + \frac{\psi^{n+1} + \psi^n}{2} = \Delta \mu^{n+1}, \tag{3.14}$$

$$\mu^{n+1} = -\epsilon^2 \Delta \frac{\phi^{n+1} + \phi^n}{2} + \frac{U^{n+1} + U^n}{2} H^\star + \beta \frac{\phi^{n+1} - \phi^n}{\delta t}, \tag{3.15}$$

$$U^{n+1} - U^n = \frac{1}{2} H^\star (\phi^{n+1} - \phi^n), \tag{3.16}$$

$$\frac{\psi^{n+1} + \psi^n}{2} = \frac{\phi^{n+1} - \phi^n}{\delta t}, \tag{3.17}$$

where

$$H^\star = \frac{f(\phi^\star)}{\sqrt{F(\phi^\star) + B}}, \quad \phi^\star = \frac{3}{2} \phi^n - \frac{1}{2} \phi^{n-1}. \tag{3.18}$$

The boundary conditions are either

(3.19) $\quad$ (i) $\phi^{n+1}, \mu^{n+1}$ are periodic; or (ii) $\partial_{\mathbf{n}} \phi^{n+1}|_{\partial\Omega} = \partial_{\mathbf{n}} \mu^{n+1}|_{\partial\Omega} = 0$.

Since the nonlinear coefficient $H$ of the new variables $U$ are treated explicitly, we can rewrite the equations (3.16) and (3.17) as follows:

$$\begin{cases} U^{n+1} = \dfrac{H^\star}{2} \phi^{n+1} + g_1^n, \\ \psi^{n+1} = \dfrac{2}{\delta t} \phi^{n+1} + g_2^n, \end{cases} \tag{3.20}$$

where $g_1^n = (U^n - \frac{H^\star}{2} \phi^n)$, $g_2^n = (-\frac{2}{\delta t} \phi^n - \psi^n)$. Thus (3.14)-(3.15) can be rewritten as the following linear system

$$\widehat{\alpha} \phi^{n+1} = \Delta \mu^{n+1} + g_3^n, \tag{3.21}$$
$$\mu^{n+1} = P_1(\phi^{n+1}) + g_4^n, \tag{3.22}$$

where

$$\begin{cases} \widehat{\alpha} = (\dfrac{\alpha}{\delta t} + \dfrac{1}{2}) \dfrac{2}{\delta t}, \\ P_1(\phi^{n+1}) = -\dfrac{\epsilon^2}{2} \Delta \phi^{n+1} + \dfrac{1}{4} H^\star H^\star \phi^{n+1} + \dfrac{\beta}{\delta t} \phi^{n+1}, \\ g_3^n = -(\dfrac{\alpha}{\delta t} + \dfrac{1}{2}) g_2^n + (\dfrac{\alpha}{\delta t} - \dfrac{1}{2}) \psi^n, \\ g_4^n = \dfrac{-\epsilon^2}{2} \Delta \phi^n + \dfrac{1}{2} H^\star (g_1^n + U^n) - \dfrac{\beta}{\delta t} \phi^n. \end{cases} \tag{3.23}$$

Therefore, we can solve $\phi^{n+1}$ and $\mu^{n+1}$ directly from (3.21) and (3.22). Once we obtain $\phi^{n+1}$, the $\psi^{n+1}, U^{n+1}$ are automatically given in (3.20). Furthermore, we notice

$$(P_1(\phi), \psi) = \frac{\epsilon^2}{2} (\nabla \phi, \nabla \psi) + \frac{1}{4} (H^\star \phi, H^\star \psi) + \frac{\beta}{\delta t} (\phi, \psi), \tag{3.24}$$

if $\psi$ enjoys the same boundary condition as $\phi$ in (3.19). Therefore, the linear operator $P_1(\phi)$ is symmetric (self-adjoint). Moreover, for any $\phi$ with $\int_\Omega \phi d\boldsymbol{x} = 0$, we have

$$(P_1(\phi), \phi) = \frac{\epsilon^2}{2} \|\nabla \phi\|^2 + \frac{1}{4} \|H^n \phi\|^2 + \frac{\beta}{\delta t} \|\phi\|^2 \geq 0, \tag{3.25}$$

where " = " is valid if and only if $\phi \equiv 0$.

We first show the well-posedness of the linear system (3.14)-(3.17) (or (3.21)-(3.22)) as follows.

**Theorem 3.1.** *The linear system* (3.21)-(3.22) *admits a unique solution in* $H^1(\Omega)$, *and the linear operator is symmetric positive definite.*

*Proof.* From (3.14), by taking the $L^2$ inner product with 1 and notice $\psi^0 = 0$, we derive

$$(3.26) \qquad (\frac{\alpha}{\delta t} + \frac{1}{2}) \int_\Omega \psi^{n+1} d\boldsymbol{x} = (\frac{\alpha}{\delta t} - \frac{1}{2}) \int_\Omega \psi^n d\boldsymbol{x} = 0.$$

From (3.17), we have

$$(3.27) \qquad \int_\Omega \phi^{n+1} d\boldsymbol{x} = \int_\Omega \phi^n d\boldsymbol{x} = \cdots = \int_\Omega \phi^0 d\boldsymbol{x}.$$

Let $V_\phi = \frac{1}{|\Omega|} \int_\Omega \phi^0 d\boldsymbol{x}$, $V_\mu = \frac{1}{|\Omega|} \int_\Omega \mu^{n+1} d\boldsymbol{x}$, and we define

$$(3.28) \qquad \widehat{\phi}^{n+1} = \phi^{n+1} - V_\phi, \widehat{\mu}^{n+1} = \mu^{n+1} - V_\mu.$$

Thus, from (3.21)-(3.22), $(\widehat{\phi}^{n+1}, \widehat{\mu}^{n+1})$ are the solutions for the following equations with unknowns $(\phi, w)$,

$$(3.29) \qquad \widehat{\alpha}\phi - \Delta w = f^n,$$
$$(3.30) \qquad w + V_\mu - P_1(\phi) = g^n,$$

where $f^n = g_3^n - \widehat{\alpha}V_\phi$, $\int_\Omega f^n d\boldsymbol{x} = 0$, $g^n = g_4^n + \frac{1}{4}H^\star H^\star \alpha_0 + \frac{\beta}{\delta t}\alpha_0$, $\int_\Omega \phi d\boldsymbol{x} = 0$ and $\int_\Omega w d\boldsymbol{x} = 0$.

Applying $-\Delta^{-1}$ to (3.29) and using (3.30), we obtain

$$(3.31) \qquad -\widehat{\alpha}\Delta^{-1}\phi + P_1(\phi) - V_\mu = -\Delta^{-1} f^n - g^n.$$

Let us express the above linear system (3.31) as $\mathbb{A}\phi = b$,

(i). For any $\phi_1$ and $\phi_2$ in $H^1(\Omega)$ satisfy the boundary conditions (2.9) and $\int_\Omega \phi_1 dx = \int_\Omega \phi_2 dx = 0$, using integration by parts, we derive

$$(3.32) \qquad \begin{aligned} (\mathbb{A}(\phi_1), \phi_2) &= -\widehat{\alpha}(\Delta^{-1}\phi_1, \phi_2) + (P_1(\phi_1), \phi_2) \\ &\leq C_1(\|\nabla\Delta^{-1}\phi_1\|\|\nabla\Delta^{-1}\phi_2\| + \|\nabla\phi_1\|\|\nabla\phi_2\| + \|\phi_1\|\|\phi_2\|) \\ &\leq C_2 \|\phi_1\|_{H^1} \|\phi_2\|_{H^1}. \end{aligned}$$

Therefore, the bilinear form $(\mathbb{A}(\phi_1), \phi_2)$ is bounded $\forall \phi_1, \phi_2 \in H^1(\Omega)$.

(ii). For any $\phi \in H^1(\Omega)$, it is easy to derive that, ,

$$(3.33) \qquad (\mathbb{A}(\phi), \phi) = \widehat{\alpha}\|\nabla\Delta^{-1}\phi\|^2 + \frac{\epsilon^2}{2}\|\nabla\phi\|^2 + \frac{1}{4}\|H^\star \phi\|^2 + \frac{\beta}{\delta t}\|\phi\|^2 \geq C_3 \|\phi\|_{H^1}^2,$$

for $\int_\Omega \phi dx = 0$ from Poincare inequality. Thus the bilinear form $(\mathbb{A}(\phi), \psi)$ is coercive.

Then from the Lax-Milgram theorem, we conclude the linear system (3.31) admits a unique solution in $H^1(\Omega)$.

For any $\phi_1, \phi_2$ with $\int_\Omega \phi_1 d\boldsymbol{x} = 0$ and $\int_\Omega \phi_2 d\boldsymbol{x} = 0$, we can easily derive

$$(3.34) \qquad (\mathbb{A}\phi_1, \phi_2) = (\phi_1, \mathbb{A}\phi_2).$$

Thus $\mathbb{A}$ is self-adjoint. Meanwhile, from (3.33), we derive $(\mathbb{A}\phi, \phi) \geq 0$, where "=" is valid if only if $\phi = 0$. This concludes the linear operator $\mathbb{A}$ is positive definite. $\square$

The energy stability of the scheme (3.14)-(3.17) (or (3.21)-(3.22)) is presented as follows.

**Theorem 3.2.** *The scheme (3.14)-(3.17) (or (3.21)-(3.22)) is unconditionally energy stable satisfying the following discrete energy dissipation law,*

$$\text{(3.35)} \quad \frac{1}{\delta t}(E_{cn2}^{n+1} - E_{cn2}^n) = -\left\|\frac{\nabla(p^{n+1} + p^n)}{2}\right\|^2 - \beta\left\|\frac{\phi^{n+1} - \phi^n}{\delta t}\right\|^2 \leq 0,$$

*where*

$$\text{(3.36)} \quad E_{cn2} = \frac{\epsilon^2}{2}\|\nabla\phi\|^2 + \|U\|^2 + \frac{\alpha}{2}\|\nabla p\|^2 - B|\Omega|.$$

*Proof.* First, we combine (3.14) and (3.15) together and apply the $\Delta^{-1}$ to obtain

$$\text{(3.37)} \quad \frac{\alpha}{\delta t}\Delta^{-1}(\psi^{n+1} - \psi^n) + \Delta^{-1}\frac{\psi^{n+1} + \psi^n}{2} = -\epsilon^2\Delta\frac{\phi^{n+1} + \phi^n}{2} + \frac{U^{n+1} + U^n}{2}H^\star + \beta\frac{\phi^{n+1} - \phi^n}{\delta t}.$$

Secondly, by taking the $L^2$ inner product of (3.37) with $\phi^{n+1} - \phi^n$, we obtain

$$\text{(3.38)} \quad \frac{\alpha}{\delta t}(\Delta^{-1}(\psi^{n+1} - \psi^n), \phi^{n+1} - \phi^n) + \frac{1}{2}(\Delta^{-1}(\psi^{n+1} + \psi^n), \phi^{n+1} - \phi^n)$$
$$= \frac{\epsilon^2}{2}(\|\nabla\phi^{n+1}\|^2 - \|\phi^n\|^2) + (\frac{U^{n+1} + U^n}{2}H^\star, \phi^{n+1} - \phi^n) + \frac{\beta}{\delta t}\|\phi^{n+1} - \phi^n\|^2.$$

Thirdly, by taking the $L^2$ inner product of (3.16) with $-(U^{n+1} + U^n)$, we obtain

$$\text{(3.39)} \quad -(\|U^{n+1}\|^2 - \|U^n\|^2) = -(\frac{1}{2}H^\star(\phi^{n+1} - \phi^n), U^{n+1} + U^n).$$

Fourthly, define $p^{n+1} = \Delta^{-1}\psi^{n+1}$. By subtracting with the $n$-step, we obtain

$$\text{(3.40)} \quad \Delta(p^{n+1} - p^n) = \psi^{n+1} - \psi^n.$$

From (3.17) and (3.40), we derive

$$\text{(3.41)} \quad \begin{aligned}\frac{\alpha}{\delta t}(\Delta^{-1}(\psi^{n+1} - \psi^n), \phi^{n+1} - \phi^n) &= \frac{\alpha}{2}(p^{n+1} - p^n, \psi^{n+1} + \psi^n) \\ &= \frac{\alpha}{2}(p^{n+1} - p^n, \Delta(p^{n+1} + p^n)) \\ &= -(\frac{\alpha}{2}\|\nabla p^{n+1}\|^2 - \frac{\alpha}{2}\|\nabla p^n\|^2),\end{aligned}$$

and

$$\text{(3.42)} \quad \begin{aligned}\frac{1}{2}(\Delta^{-1}(\psi^{n+1} + \psi^n), \phi^{n+1} - \phi^n) &= \frac{\delta t}{4}(p^{n+1} + p^n, \psi^{n+1} + \psi^n) \\ &= \frac{\delta t}{4}(p^{n+1} + p^n, \Delta(p^{n+1} + p^n)) \\ &= -\frac{\delta t}{4}\|\nabla(p^{n+1} + p^n)\|^2.\end{aligned}$$

Finally, by combining (3.38), (3.39), (3.41) and (3.42), we obtain

$$\text{(3.43)} \quad \begin{aligned}\frac{\epsilon^2}{2}(\|\nabla\phi^{n+1}\|^2 - \|\nabla\phi^n\|^2) + \|U^{n+1}\|^2 - \|U^n\|^2 + \frac{\alpha}{2}(\|\nabla p^{n+1}\|^2 - \|\nabla p^n\|^2) \\ = -\frac{\delta t}{4}\|\nabla(p^{n+1} + p^n)\|^2 - \frac{\beta}{\delta t}\|\phi^{n+1} - \phi^n\|^2,\end{aligned}$$

that concludes the theorem. □

**Remark 3.4.** *The proposed scheme (3.14)-(3.17) follows the new energy dissipation law (3.12) formally instead of the energy law for the originated system (2.15). In the time-discrete case, the energy $E(\phi^{n+1}, U^{n+1})$ (defined in (3.36)) can be rewritten as a second order approximation to the Lyapunov functionals in $E(\phi^{n+1})$ (defined in (2.15)), that can be observed from the following facts, heuristically. Assuming the case for double well potential, from (3.16), we have*

$$(3.44) \qquad U^{n+1} - (\sqrt{F(\phi^{n+1}) + B}) = U^n - (\sqrt{F(\phi^n) + B}) + R_{n+1},$$

*where $R_{n+1} = O((\phi^{n+1} - \phi^n)(\phi^{n+1} - 2\phi^n + \phi^{n-1}))$. Since $R_k = O(\delta t^3)$ for $0 \leq k \leq n+1$ and $U^0 = (\sqrt{F(\phi^0) + B})$, by mathematical induction we can easily get*

$$(3.45) \qquad U^{n+1} = \sqrt{F(\phi^{n+1}) + B} + O(\delta t^2).$$

3.2. **Adam-Bashforth Scheme.** Next, for the completion of the development of second order schemes, we further develop a scheme based on the Adam-Bashforth backward differentiation formula (BDF2). It provides an alternative second order scheme with the unconditional energy stability that is beneficial for the scheme development. Since the stability proof of the BDF2 scheme is quite different from the Crank-Nicolson scheme, we give its details as well.

<u>Scheme</u> **2.** *Given the initial condition ($\phi^0$, $U^0$), compute $U^1$ and $\phi^1$ by assuming $U^{-1} = U^0$ and $\phi^{-1} = \phi^0$ for the initial step. Having computed ($\phi^n$, $U^n$) and ($\phi^{n-1}$, $U^{n-1}$), with $n \geq 1$, we solve $\phi^{n+1}$, $U^{n+1}$ as follows:*

$$(3.46) \qquad \alpha \frac{3\psi^{n+1} - 4\psi^n + \psi^{n-1}}{2\delta t} + \psi^{n+1} = \Delta \mu^{n+1},$$

$$(3.47) \qquad \mu^{n+1} = -\epsilon^2 \Delta \phi^{n+1} + U^{n+1} H^\dagger + \beta \frac{3\phi^{n+1} - 4\phi^n + \phi^{n-1}}{2\delta t},$$

$$(3.48) \qquad 3U^{n+1} - 4U^n + U^{n-1} = \frac{1}{2} H^\dagger (3\phi^{n+1} - 4\phi^n + \phi^{n-1}),$$

$$(3.49) \qquad \psi^{n+1} = \frac{3\phi^{n+1} - 4\phi^n + \phi^{n-1}}{2\delta t},$$

*where*

$$(3.50) \qquad H^\dagger = \frac{f(\phi^\dagger)}{\sqrt{F(\phi^\dagger) + B}}, \quad \phi^\dagger = 2\phi^n - \phi^{n-1}.$$

*The boundary conditions are*

$$(3.51) \qquad (i) \ \phi^{n+1}, \mu^{n+1} \text{ are periodic; or } (ii) \ \partial_{\mathbf{n}} \phi^{n+1}|_{\partial \Omega} = \partial_{\mathbf{n}} \mu^{n+1}|_{\partial \Omega} = 0.$$

Similar to the Crank-Nicolson scheme, we can rewrite the equations (3.48) and (3.49) as follows:

$$(3.52) \qquad \begin{cases} U^{n+1} = \dfrac{H^\dagger}{2} \phi^{n+1} + h_1^n, \\ \psi^{n+1} = \dfrac{3}{2\delta t} \phi^{n+1} + h_2^n, \end{cases}$$

where $h_1^n = (U^\pm - \frac{H^\dagger}{2}\phi^\pm)$, $h_2^n = \frac{3}{2\delta t}\phi^\pm$ with $S^\pm = \frac{4S^n - S^{n-1}}{3}$ for any variable $S$. Thus (3.46)-(3.47) can be rewritten as the following linear system

$$(3.53) \qquad \widetilde{\alpha} \phi^{n+1} = \Delta \mu^{n+1} + h_3^n,$$

$$(3.54) \qquad \mu^{n+1} = P_2(\phi^{n+1}) + h_4^n,$$

where

(3.55)
$$\begin{cases} P_2(\phi^{n+1}) = -\epsilon^2 \Delta \phi^{n+1} + \frac{1}{2}H^\dagger H^\dagger \phi^{n+1} + \frac{3\beta}{2\delta t}\phi^{n+1}, \\ h_3^n = -(\frac{3\alpha}{2\delta t} + 1)h_2^n + \frac{3\alpha}{2\delta t}\psi^\pm, \\ h_4^n = \frac{1}{2}H^\dagger h_1^n - \frac{3\beta}{2\delta t}\phi^\pm, \\ \widetilde{\alpha} = (\frac{3\alpha}{2\delta t} + 1)\frac{3}{2\delta t}. \end{cases}$$

Actually, we can solve $\phi^{n+1}$ and $\mu^{n+1}$ directly from (3.53) and (3.54). Once we obtain $\phi^{n+1}$, the $\psi^{n+1}, U^{n+1}$ is automatically given in (3.52). Furthermore, we notice

(3.56) $$(P_2(\phi), \psi) = \epsilon^2(\nabla\phi, \nabla\psi) + \frac{1}{2}(H^\dagger\phi, H^\dagger\psi) + \frac{3\beta}{2\delta t}(\phi, \psi),$$

if $\psi$ enjoys the same boundary condition as $\phi$ in (3.51). Therefore, the linear operator $P_2(\phi)$ is symmetric (self-adjoint). Moreover, for any $\phi$ with $\int_\Omega \phi d\boldsymbol{x} = 0$, we have

(3.57) $$(P_2(\phi), \phi) = \epsilon^2 \|\nabla\phi\|^2 + \frac{1}{2}\|H^\dagger\phi\|^2 \geq 0,$$

where " $=$ " is valid if and only if $\phi \equiv 0$.

**Remark 3.5.** *One can show the well-posedness of the linear system (3.46)-(3.49) (or (3.53)-(3.54)). Likewise, when we rewrite (3.53)-(3.54) into a linear equation using the inverse Laplace operator, we can show the linear operator is symmetric (self-adjoint) and positive definite.*

The energy stability of the scheme (3.46)-(3.49) (or (3.53)-(3.54)) is presented as follows.

**Theorem 3.3.** *The scheme (3.46)-(3.49) (or (3.53)-(3.54)) is unconditionally energy stable satisfying the following discrete energy dissipation law,*

(3.58) $$\frac{1}{\delta t}(E_{bdf2}^{n+1} - E_{bdf2}^n) \leq -\|\nabla p^{n+1}\|^2 - \beta\left\|\frac{3\phi^{n+1} - 4\phi^n + \phi^{n-1}}{2\delta t}\right\|^2 \leq 0,$$

where

(3.59) $$E_{bdf2}^{n+1} = \frac{\epsilon^2}{2}\Big(\frac{\|\nabla\phi^{n+1}\|^2 + \|2\nabla\phi^{n+1} - \nabla\phi^n\|^2}{2}\Big) + \frac{\|U^{n+1}\|^2 + \|2U^{n+1} - U^n\|^2}{2} \\ + \frac{\alpha}{2}\frac{\|\nabla p^{n+1}\|^2 + \|2\nabla p^{n+1} - \nabla p^n\|^2}{2} - B|\Omega|.$$

*Proof.* First, from (3.46), by taking the $L^2$ inner product with 1 and notice $\psi^0 = 0$, we derive

(3.60) $$(\frac{3\alpha}{2\delta t} + 1)\int_\Omega \psi^{n+1} d\boldsymbol{x} = \frac{4\alpha}{2\delta t}\int_\Omega \psi^n d\boldsymbol{x} - \frac{\alpha}{2\delta t}\int_\Omega \psi^{n-1} d\boldsymbol{x} = 0.$$

where we use $\int_\Omega \psi^1 d\boldsymbol{x} = 0$, this is valid since $\psi^1$ can be obtained using the Crank-Nicolson scheme.

Second, we combine (3.46) and (3.47) together and applying the $\Delta^{-1}$ to obtain

(3.61) $$\alpha \Delta^{-1}(\frac{3\psi^{n+1} - 4\psi^n + \psi^{n-1}}{2\delta t}) + \Delta^{-1}\psi^{n+1} \\ = -\epsilon^2 \Delta\phi^{n+1} + U^{n+1}H^\dagger + \beta\frac{3\phi^{n+1} - 4\phi^n + \phi^{n-1}}{2\delta t}.$$

Third, by taking the $L^2$ inner product of (3.37) with $3\phi^{n+1} - 4\phi^n + \phi^{n-1}$, and applying the following identity

(3.62) $$(3a - 4b + c, 2a) = a^2 - b^2 + (2a - b)^2 - (2b - c)^2 + (a - 2b + c)^2,$$

we obtain

(3.63) $$\begin{aligned}\frac{\alpha}{2\delta t}(\Delta^{-1}(3\psi^{n+1} - 4\psi^n + \psi^{n-1}), 3\phi^{n+1} - 4\phi^n + \phi^{n-1}) + (\Delta^{-1}\psi^{n+1}, 3\phi^{n+1} - 4\phi^n + \phi^{n-1}) \\ = \frac{\epsilon^2}{2}\Big(\|\nabla\phi^{n+1}\|^2 - \|\nabla\phi^n\|^2 + \|2\nabla\phi^{n+1} - \nabla\phi^n\|^2 - \|2\nabla\phi^n - \nabla\phi^{n-1}\|^2 \\ + \|3\nabla\phi^{n+1} - 4\nabla\phi^n - \nabla\phi^{n-1}\|^2\Big) \\ + (U^{n+1}H^\dagger, 3\phi^{n+1} - 4\phi^n + \phi^{n-1}) + \frac{\beta}{2\delta t}\|3\phi^{n+1} - 4\phi^n + \phi^{n-1}\|^2.\end{aligned}$$

Third, by taking the $L^2$ inner product of (3.16) with $-2U^{n+1}$, we obtain

(3.64) $$\begin{aligned}-(\|U^{n+1}\|^2 - \|U^n\|^2 + \|2U^{n+1} - U^n\|^2 - \|2U^n - U^{n-1}\|^2 + \|U^{n+1} - 2U^n + U^{n-1}\|^2) \\ = -(H^\dagger(3\phi^{n+1} - 4\phi^n + \phi^{n-1}), U^{n+1}).\end{aligned}$$

Fourth, define $p^{n+1} = \Delta^{-1}\psi^{n+1}$, by subtracting with the $n$ and $n-1$-step, we obtain

(3.65) $$\Delta(3p^{n+1} - 4p^n + p^{n-1}) = 3\psi^{n+1} - 4\psi^n + \psi^{n-1}.$$

From (3.49) and (3.65), we derive

(3.66) $$\begin{aligned}\frac{\alpha}{2\delta t}(\Delta^{-1}(3\psi^{n+1} - 4\psi^n + \psi^{n-1}), 3\phi^{n+1} - 4\phi^n + \phi^{n-1}) \\ = \alpha(3p^{n+1} - 4p^n + p^{n-1}, \psi^{n+1}) \\ = \alpha(3p^{n+1} - 4p^n + p^{n-1}, \Delta p^{n+1}) \\ = -\frac{\alpha}{2}(\|\nabla p^{n+1}\|^2 - \|\nabla p^n\|^2 + \|2\nabla p^{n+1} - \nabla p^n\|^2 - \|2\nabla p^n - \nabla p^{n-1}\|^2 \\ + \|\nabla p^{n+1} - 2\nabla p^n + \nabla p^{n-1}\|^2),\end{aligned}$$

and

(3.67) $$\begin{aligned}(\Delta^{-1}\psi^{n+1}, 3\phi^{n+1} - 4\phi^n + \phi^{n-1}) &= 2\delta t(p^{n+1}, \psi^{n+1}) \\ &= 2\delta t(p^{n+1}, \Delta p^{n+1}) \\ &= -2\delta t\|\nabla p^{n+1}\|^2.\end{aligned}$$

Finally, by combining (3.63), (3.64), (3.66) and (3.67), we obtain

$$\frac{\epsilon^2}{2}(\|\nabla\phi^{n+1}\|^2 - \|\nabla\phi^n\|^2 + \|2\nabla\phi^{n+1} - \nabla\phi^n\|^2 - \|2\nabla\phi^n - \nabla\phi^{n-1}\|^2 + \|\nabla\phi^{n+1} - 2\nabla\phi^n + \nabla\phi^{n-1}\|^2)$$
$$+ \|U^{n+1}\|^2 - \|U^n\|^2 + \|2U^{n+1} - U^n\|^2 - \|2U^n - U^{n-1}\|^2 + \|U^{n+1} - 2U^n + U^{n-1}\|^2$$
$$+ \frac{\alpha}{2}(\|\nabla p^{n+1}\|^2 - \|\nabla p^n\|^2 + \|2\nabla p^{n+1} - \nabla p^n\|^2 - \|2\nabla p^n - \nabla p^{n-1}\|^2 + \|\nabla p^{n+1} - 2\nabla p^n + \nabla p^{n-1}\|^2)$$
$$= -2\delta t\|\nabla p^{n+1}\|^2 - \frac{\beta}{2\delta t}\|3\phi^{n+1} - 4\phi^n + \phi^{n-1}\|^2.$$

That concludes the theorem. $\square$

**Remark 3.6.** *Heuristically, the $\frac{1}{\delta t}(E_{bdf2}^{n+1} - E_{bdf2}^n)$ is a second order approximation of $\frac{d}{dt}E(\phi, U)$ at $t = t^{n+1}$. For instance, for any smooth variable $S$ with time, one can write*
$$\left(\frac{\|S^{n+1}\|^2 + \|2S^{n+1} - S^n\|^2}{2\delta t}\right) - \left(\frac{\|S^n\|^2 + \|2S^n - S^{n-1}\|^2}{2\delta t}\right)$$
$$\cong \left(\frac{\|S^{n+2}\|^2 - \|S^n\|^2}{2\delta t}\right) + O(\delta t^2) \cong \frac{d}{dt}\|S(t^{n+1})\|^2 + O(\delta t^2).$$

**Remark 3.7.** *The new variable $U$ is introduced in order to handle the nonlinear bulk potential $F(\phi)$. Since the discrete energy still includes the gradient term of $\phi$, therefore, due to the Poincare inequality and mass conservation property for Cahn-Hilliard equation $\int_\Omega \phi^{n+1} d\boldsymbol{x} = \int_\Omega \phi^0 d\boldsymbol{x}$, the $H^1$ stability for the variable $\phi$ is still valid for the proposed scheme, which makes it possible to implement the rigorous error analysis. About the complete error analysis of the IEQ type schemes for solving the classical CH equation with general nonlinear bulk potentials, we refer to a recent article [67], in which, some reasonable sufficient conditions about boundedness and continuity for the nonlinear potential are given and optimal error estimates are obtained. Similar work can be performed for the VCH-HR system as well with no essential difficulties.*

## 4. Numerical tests.

In this section, we present various numerical experiments to validate the theories derived in the previous section and demonstrate the efficiency, energy stability and accuracy of the proposed numerical schemes. In all examples, we set the domain $\Omega = [0,1]^d, d = 2, 3$ unless elaborated. We use the second order central finite difference method to discretize the space operators in the semi-discretized model. In all simulations, we set $\epsilon = 0.01$, and $\alpha, \beta$ will be chosen accordingly. For double well potential case, we set $B = 0$ and $U = \phi(1-\phi)$. For Flory-Huggins case, we set $B = 1$, $\sigma = 0.001$ and $\chi = 2.5$.

### 4.1. Convergence test.
We first test the convergence rates of the two proposed schemes, the second-order Crank-Nicolson scheme (CN2) (3.14)-(3.17) and the second-order Backward-Difference scheme (BDF2) (3.46)-(3.49). Use the following initial condition

$$\phi = 0.5\Big(1 + \max(\tanh\frac{0.2 - R_1}{\varepsilon}, \tanh\frac{0.2 - R_2}{\varepsilon})\Big), \tag{4.1}$$

with $R_1 = \sqrt{(x-0.71)^2 + (y-0.5)^2}$ and $R_2 = \sqrt{(x-0.29)^2 + (y-0.5)^2}$. The initial profile of $\phi$ is shown in the first panel of Fig. 2. The spacial mesh is $256 \times 256$. We perform the time-step refinement test to obtain the order of convergence in time by taking a linear refinement path for time step $\delta t = \frac{0.01}{2^k}, k = 0, 1, \cdots, 6$. The numerical errors are calculated as the difference between the solution of coarse time step and that of the adjacent finner time step. We plot the Cauchy sequence of $L^2$ errors at $t = 4$ with different time step sizes in Fig. 1 and the convergence rate is shown to be second order for both schemes.

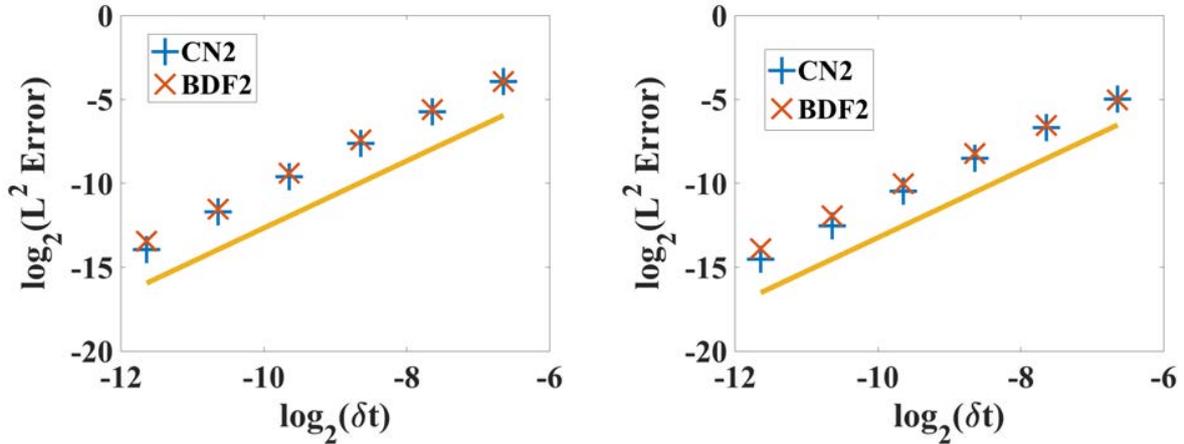

(a) $L^2$ error for the phase variable $\phi$.

(b) $L^2$ error for the auxillary variable $U$.

FIG. 1. Convergence test for the $L^2$ errors for $\phi$ and $U$ computed by the the second order scheme CN2 and BDF2 using different temporal resolutions at $t = 4$. The time step is $\delta t = 0.01(\frac{1}{2})^k$ for $k = 0, 1, 2, 3, 4, 5, 6$ and the numerical errors are calculated as the difference between the solution of the coarse time step and that of the adjacent finner time step.

4.2. **The viscous and hyperbolic relaxation effects for the coalescence of two kissing bubbles.** In this example, we study the coalescence dynamics of two kissing bubbles by varying the viscous and hyperbolic relaxation parameters $\alpha$ and $\beta$. The computational domain $\Omega$ is still $[0, 1]^2$ and the initial profile is given in (4.1), and we use the CN2 scheme and $128^2$ grid points to discretize the domain.

We start with the classical Cahn-Hilliard equation by setting $\alpha = \beta = 0$. In Fig. 2, the two bubbles coalesces into one big bubble ( in a lower free energy state) due to the combination of the surface tension effect. Then, we further set the hyperbolic relaxation parameter $\alpha = 1$ while keeping viscous parameter of $\beta = 0$. The numerical result in shown in Fig. 3. At $t = 1$ and 2, the interface of the circle shows some sawtooth profile and eventually forms a circle, i.e. the Cahn-Hilliard equation with hyperbolic relaxation term predicts different dynamics, but the same final steady state. Then we investigate the viscous effect by setting $\beta = 1$ and $\alpha = 0$. The numerical results are illustrated in Fig. 4. We observe that the coalesce speed is much slower than the two cases where $\beta = 0$.

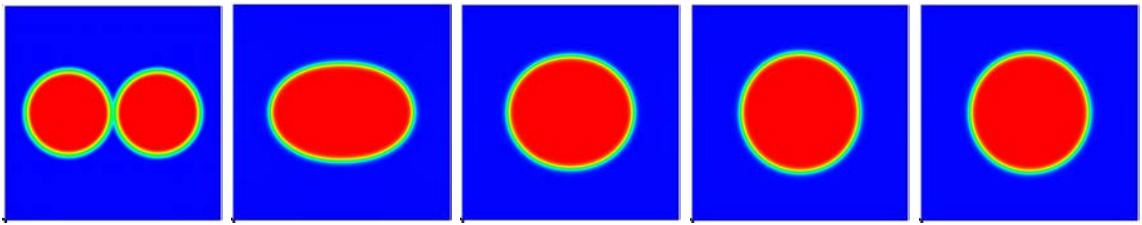

FIG. 2. Time evolution of the drop in 2D when $\beta = 0$ and $\alpha = 0$. Snapshots are taken at $t = 0, 1, 2, 10, 200$.

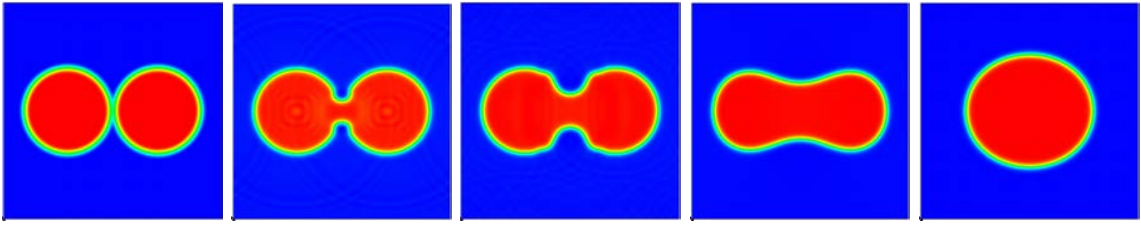

Fig. 3. Time evolution of the drop in 2D when $\beta = 0$ and $\alpha = 1$. Snapshots are taken at $t = 0, 1, 2, 10, 200$.

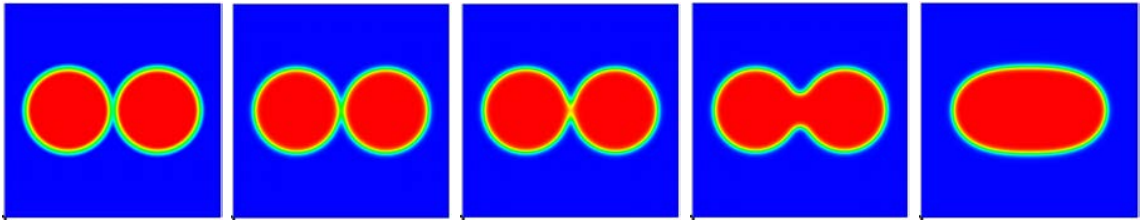

Fig. 4. Time evolution of the drop in 2D when $\beta = 1$ and $\alpha = 0$. Snapshots are taken at $t = 0, 1, 2, 10, 200$.

We plot the evolution of energy curves for nine cases in Fig. 5 when both $\alpha$ and $\beta$ take the three values of $0, 0.5, 1$. We find that $\beta$ can dramatically affect the speed of coalesces than $\alpha$. Both the results from Scheme 1 and Scheme 2 are shown, and they predict the same dynamics.

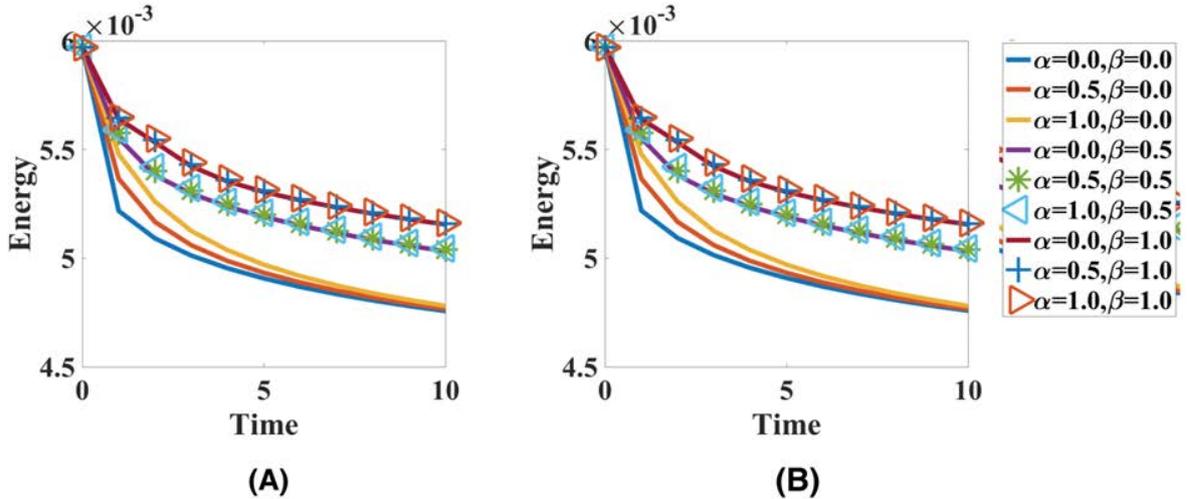

Fig. 5. Time evolution of the free energy functional for the coalescence of two kissing bubbles for nine choices of order parameters $\alpha = 0, 0.5, 1$, and $\beta = 0, 0.5, 1$. (A) energy plot of (3.11) using Scheme 1; (B) energy plot of (3.11) using Scheme 2.

4.3. **Energy.** Since the two developed schemes follow a modified energy law with the new energy (3.11) instead of the original energy (2.14). In the continuous level, these are energies are equivalent. However, when they are discretized, the new energy (3.11) is a second-order approximation of the original one (2.14), as stated in Remark 3.4 and 3.6. To verify this statement, we perform the following simulations. We use the initial condition $\phi = 0.5(1 + \cos(2\pi z)\cos(2\pi y))$ and periodic boundary conditions. Without the loss of generality, we choose $\alpha = \beta = 0.5$. Since we do not have the analytical solution, the original energy (2.14) is calculated by using the fully implicit scheme with the time step $\delta t = 2^{-12}$ as the benchmark solution. In Fig. 6, the second-order convergence of the transformed energy (3.11) that are computed by the Scheme 1 to the original energy (2.14) is observed via the time-step refinement test. To further demonstrate the effectiveness of our proposed unconditionally energy stable schemes, we further plot the time evolution of the modified energy (3.11) using progressively larger time steps with $\delta t = 0.00625, 0.01235, 0.25, 0.5$ and 1. Fig. 7 demonstrates the energy decays even with very larger time steps like $\delta t = 1$, which means the schemes behave truly as the theory indicates, i.e., they are stable for progressively large time steps.

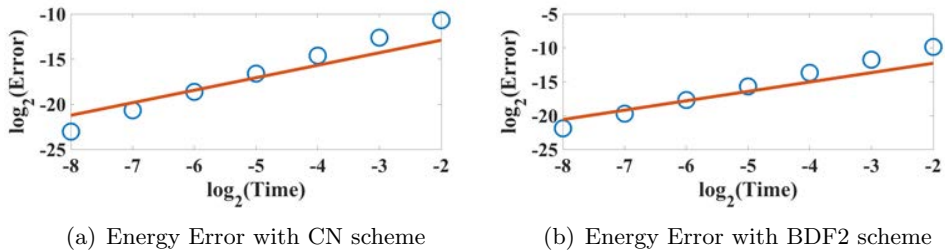

(a) Energy Error with CN scheme    (b) Energy Error with BDF2 scheme

FIG. 6. The energy difference between the original energy formulation (2.14) and the approximated energy (3.11) that is computed by Scheme 1 with different time steps. This figure shows the approximated energy (3.11) is definitely a second-order approximation of the original energy (2.14).

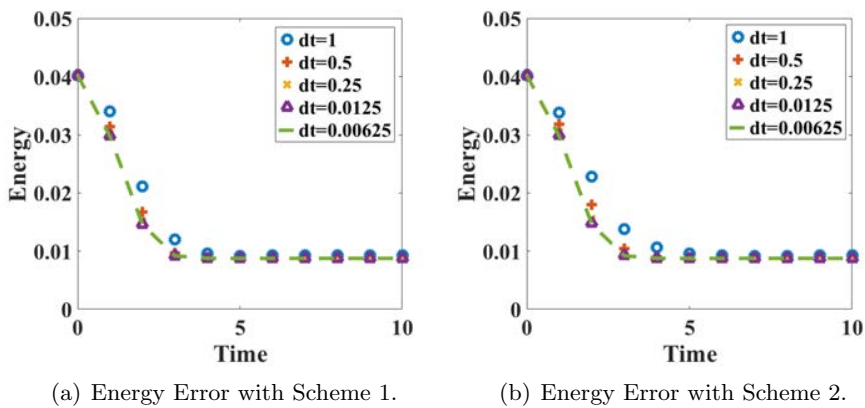

(a) Energy Error with Scheme 1.    (b) Energy Error with Scheme 2.

FIG. 7. The time evolution of the energy (3.11) using various time step sizes $\delta t = 0.00625, 0.01235, 0.25, 0.5$ and 1.

4.4. **Spinodal Decomposition in 3D.** In this example, we study the phase separation dynamics in 3D that is called " spinodal decomposition ". The process of the phase separation can be studied by considering a homogeneous binary mixture, which is quenched into the unstable part of its miscibility gap. In this case, the spinodal decomposition takes place, which manifests in the spontaneous growth of the concentration fluctuations that leads the system from the homogeneous to the two-phase state. Shortly after the phase separation starts, the domains of the binary components are formed and the interface between the two phases can be specified [1, 15, 79].

The initial conditions are taken as the randomly perturbed concentration fields as follows,

$$\phi_0(x,y,z) = \overline{\phi}_0 + 0.001\text{rand}(x,y,z), \tag{4.2}$$

where the rand$(x,y)$ represents the random number in $[0,1]$ and has zero mean. The computational domain is $[0, 2\pi]^3$ and we use the scheme CN2 and $128^3$ grid points to discretize the domain, the time step is $\delta t = 0.001$ for all 3D simulations.

From the 2D tests, we know the viscous parameter $\beta$ can have more effects on the dynamics than the hyperbolic parameter $\alpha$. Thus in the following 3D simulations, we simply set $\beta = 0.9$ and $\alpha = 0$. The red domain, corresponding to the larger values of $\phi = 1$, indicates the concentrated polymer segments [24], and the blue region, corresponding to the smaller values of $\phi = -1$, indicates the macromolecular microspheres (MMs). In Fig. 8, we perform numerical simulations for the initial profile $\overline{\phi}_0 = 0.5$, that means the volume fraction of the polymer segments are almost same as the surrounding MMs. The final steady state forms the uniform two layer structure around $t = 1800$. Fig. 9 shows the dynamical behaviors of the phase separation for the initial value $\overline{\phi}_0 = 0.3$ which means the volume of the MMS are much more than that of the polymer segments. We observe that the MMS finally accumulate together to the cylindrical shape.

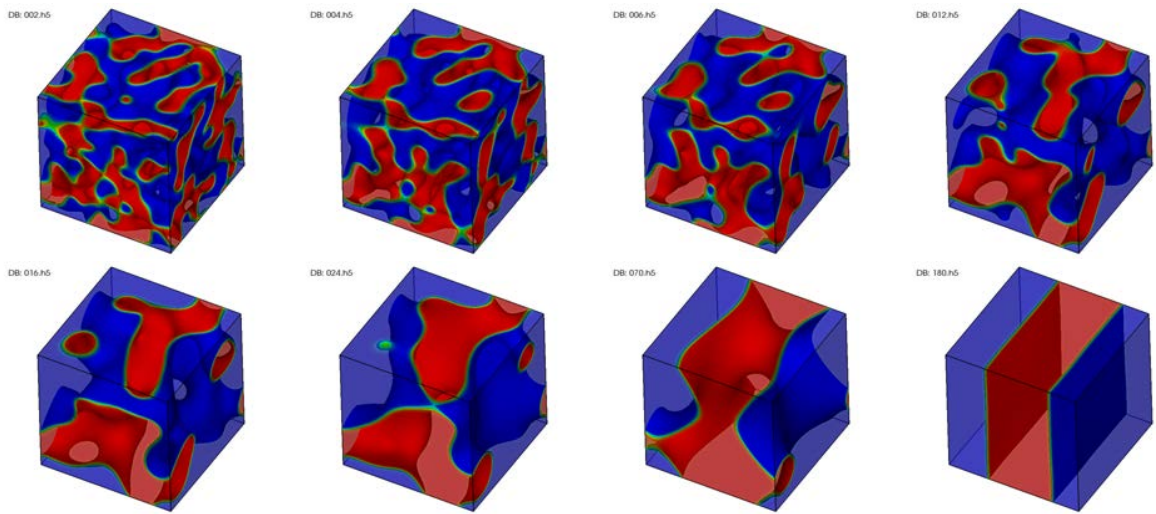

FIG. 8. 3D spinodal decomposition for random initial data with $\overline{\phi}_0 = 0.5$. Snapshots of the phase variables $\phi$ are taken at $t = 20, 40, 60, 120, 160, 240, 700,$ and $1800$. The order parameter is $\alpha = 0$ and $\beta = 0.9$.

Note the boundary conditions of the governing system can be the periodic or no-flux, in Fig. 10, we perform numerical simulations for the initial profile of $\overline{\phi}_0 = 0.3$ and $\alpha = \beta = 0$ for these two boundary conditions. For both cases, we observe that the MMS finally accumulate together to the sperical shape, where the final shape is 1/8 spherical segment at a corner for no-flux condition, and a sphere for the period boundary condition (with parts at each corner of the cube). We plot the evolution of energy curves for both cases in Fig. 11. And we observe the energy of the spherical segment with no-flux boundary condition is smaller than that of the sphere in the period boundary condition.

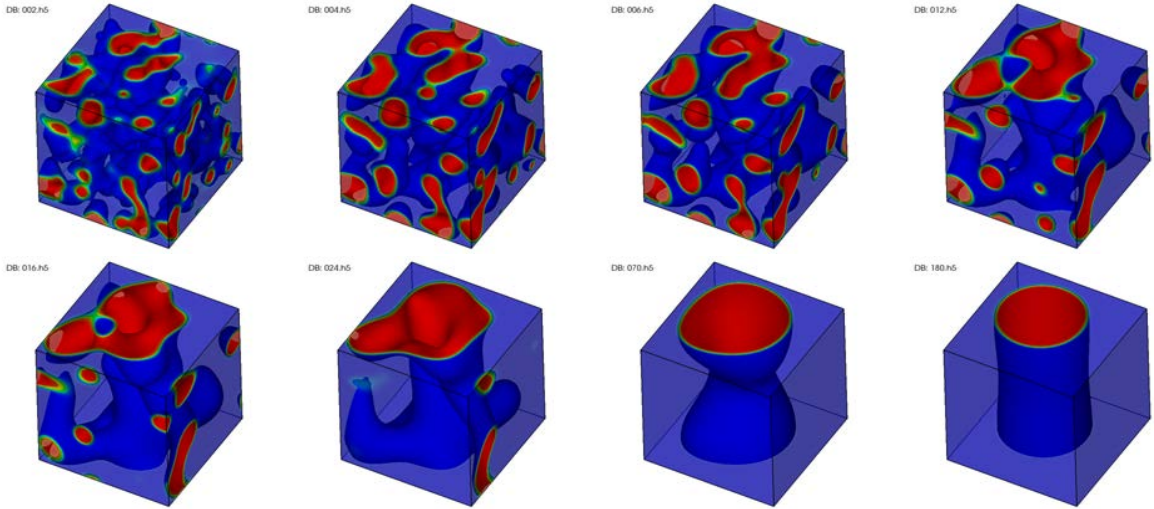

FIG. 9. 3D spinodal decomposition for random initial data with $\overline{\phi}_0 = 0.3$. Snapshots of the phase variables $\phi$ are taken at $t = 20, 40, 60, 120, 160, 240, 700$ and $2000$. The order parameter is $\alpha = 0$ and $\beta = 0.9$.

Here we also conducted several numerical tests using the Flory-Huggin potential in (3.13). The numerical results are shown in Figure 12. By using different initial value of $\overline{\phi}_0$, the predicted dynamics are dramatically different, and the steady state are distinct as well. A detailed discussion on the correlations of initial values and final steady state is out of scope for current paper. Interested readers are encouraged to conduct the numerical studies using our proposed schemes.

## 5. Concluding Remarks

In this paper, we develop two second-order in time schemes to solve the viscous Cahn-Hilliard equation with hyperbolic relaxation terms, by utilizing the novel IEQ approach. It is effective and efficient, and particularly suitable to discretize the complicated nonlinear potential with lower bound. Compared to the prevalent nonlinear schemes based on the convex splitting approaches or other nonlinear schemes, the IEQ approach can easily conquer the inconvenience from nonlinearities by linearizing the nonlinear terms in the new way. The developed schemes (i) are *accurate* (ready for second or higher order in time); (ii) are *stable* (unconditional energy dissipation law holds); and (iii) are *easy to implement* (only need to solve linear equations at each time step). Furthermore, the induced linear system is symmetric positive definite, thus one can apply any Krylov subspace methods with mass lumping as pre-conditioners for solving such system efficiently. We emphasize that, to the best of the authors' knowledge, the schemes to solve the case of logarithmic potential are the first such linear and accurate schemes with provable energy stabilities. Finally, the method is general enough to be extended to develop linear schemes for a large class of gradient flow problems with complex nonlinearities in the free energy density. Although we consider only time discrete schemes in this study, the results can be carried over to any consistent finite-dimensional Galerkin

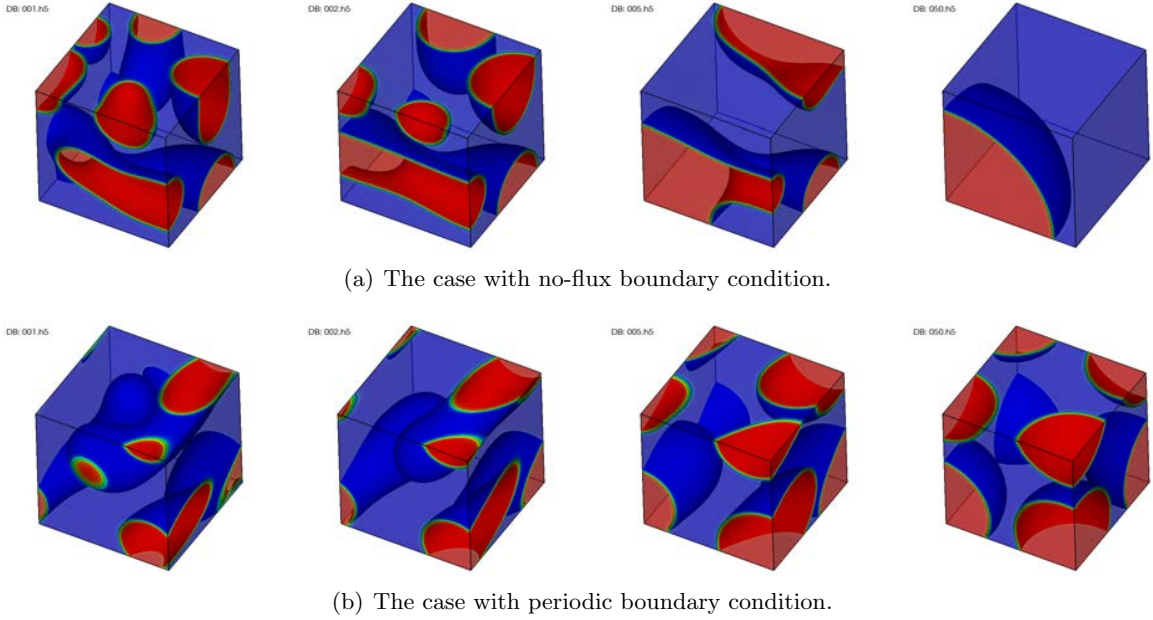

(a) The case with no-flux boundary condition.

(b) The case with periodic boundary condition.

FIG. 10. 3D spinodal decomposition for random initial data with $\overline{\phi}_0 = 0.3$ for no-flux and periodic boundary conditions. Snapshots are taken at $t = 1, 2, 5, 50$. And the order parameters are $\alpha = \beta = 0$ for both cases.

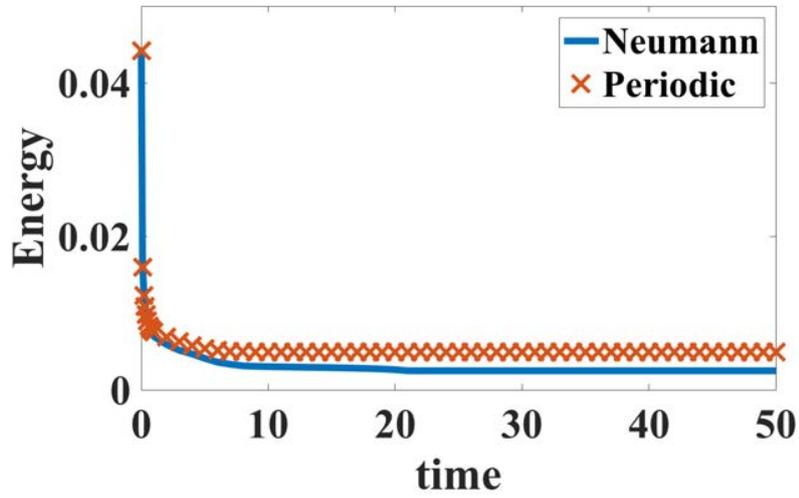

FIG. 11. Time evolution of the free energy functional for spinodal decomposition for no-flux and periodic boundary conditions with $\overline{\phi}_0 = 0.3, \alpha = \beta = 0$.

approximations since the proofs are all based on a variational formulation with all test functions in the same space as the space of the trial functions.

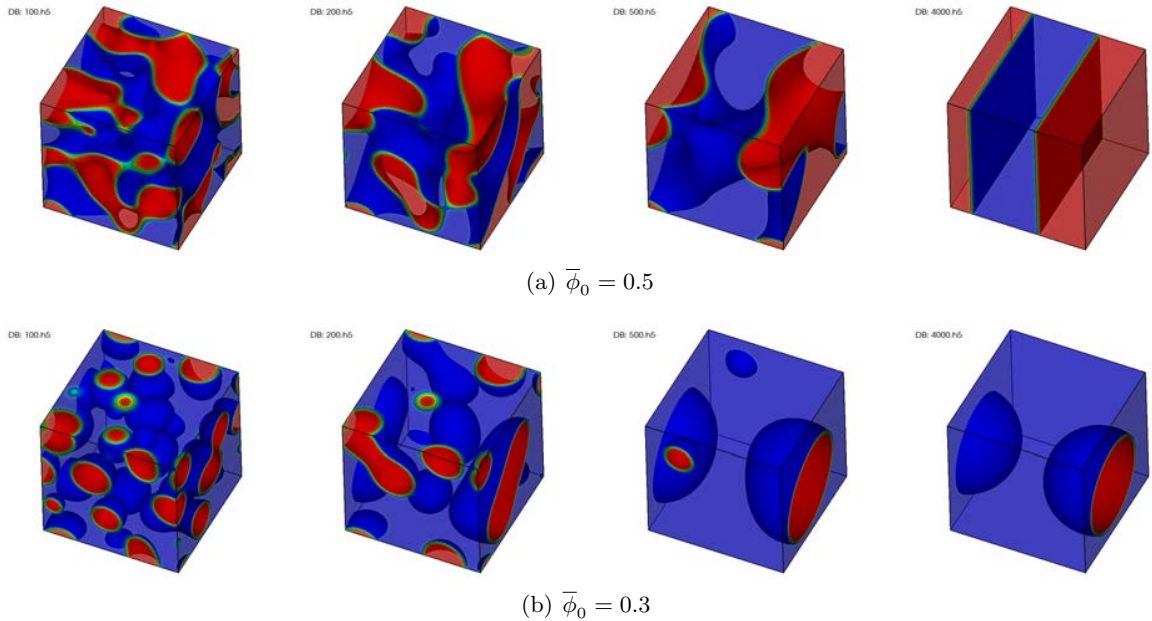

(a) $\overline{\phi}_0 = 0.5$

(b) $\overline{\phi}_0 = 0.3$

FIG. 12. 3D Spinodal decomposition with Flory-Huggins free energy. Here we choose $\alpha = \beta = 0.5$ and the initial condition as (4.2) with (a) $\overline{\phi}_0 = 0.5$ and (b) $\overline{\phi}_0 = 0.3$. The plot of $\phi$ at time $t = 100, 200, 500, 4000$ are shown respectively.

**Acknowledgments.** X. Yang is partially supported by NSF Grants DMS-1418898 and NSF-1720212. The authors thank the hospitality of Beijing Computational Science Research Center during their visits when the research was done.